\documentclass[a4paper,11pt,numreferences,mathsec,kaplist]{isueps}
\usepackage{isu}

\begin{document}
\setcounter{aqwe}{2}
\begin{article}

\begin{opening}

\udk{519.642.5}

\title{Polynomial spline collocation method for solving\\ weakly regular Volterra integral equations of\\ the first kind\thanks{This work was supported by Ministry of Science and Higher Education. The authors received no third party funding for this study. The research was carried out under State Assignment Project (no. FWEU-2021-0006) of the Fundamental Research Program of Russian Federation 2021-2030 using the resources of the High-Temperature Circuit Multi-Access Research Center (Ministry of Science and Higher Education of the Russian Federation, project no 13.CKP.21.0038). }
}

\author{A.\surname{Tynda
$^1$},  S.\surname{Noeiaghdam$^{2,3}$} and D. \surname{Sidorov$^{4,2,5}$}}

\institute{\small {$^1$ Higher and Applied Mathematics Department, Penza State University,  Penza, Russia\\
$^2$ Industrial Mathematics Laboratory, Baikal School of BRICS, Irkutsk National Research Technical University, Irkutsk, Russia\\
$^3$ Department of Applied Mathematics and Programming, South Ural State University, Chelyabinsk,  Russia\\
$^4$ Energy Systems Institute, Siberian Branch of Russian Academy of Science,  Irkutsk, Russia\\
$^5$ Institute of Mathematics and Information Technologies, Irkutsk State University,  Irkutsk, Russia
}}

\runningtitle{Weakly Regular Volterra Integral Equations of the First Kind}
\runningauthor{{A. Tynda, S. Noeiaghdam and D. Sidorov}}

\begin{abstract}
The polynomial spline collocation method is proposed for solution of  Volterra integral equations of the first kind with special piecewise continuous kernels.
The Gauss-type quadrature formula is used to approximate integrals during the discretization of the proposed projection method. The estimate of accuracy of approximate solution is
obtained. Stochastic arithmetics is also used based on the Contr\^{o}le et Estimation Stochastique des Arrondis de Calculs (CESTAC) method  and the Control of Accuracy and Debugging for Numerical Applications (CADNA) library. Applying this approach it is possible to find optimal parameters of the projective method. The numerical examples are included to illustrate the efficiency of proposed  novel collocation method.
\end{abstract}


\end{opening}
%%%%%%%%%%%%%%%%%%%%%%%%%%%%%%%%%%%%%%%%%%%%%%%%%%%%%%%%%%%%%%%%%%%%%%%%%%%%%%%%%%%%%%%%%%%

\avtogl{A. Tynda, S. Noeiaghdam, D. Sidorov} {Polynomial Spline Collocation Method for Solving Weakly Regular Volterra Integral Equations of the First Kind}

\section{Introduction}

This article focuses on the following weakly regular Volterra equations of the first kind
\begin{equation}\label{math2021-e01}
\int_{0}^{t} K(t,s) x(s) ds = g(t), \,\,\,  0 \leq s\leq t \leq T,\,\, g(0)=0,
\end{equation}
where jump discontinuous kernels are defined as follows
\begin{equation}\label{kernel}
    K(t,s) = \left\{ \begin{array}{ll}
          K_1(t,s), \,\, t,s \in m_1, \\
         \,\, \dots \,\, \dots \dots \dots  \\
         K_n(t,s), \,\, t,s \in m_n, \\
        \end{array} \right. \,\,
\end{equation}
where $m_i = \{ t, s\,\,  \bigl |\,\, \alpha_{i-1}(t) < s < \alpha_i(t) \} $,
$ { \alpha_0(t)=0,\,\, \alpha_n(t)=t,\, i=\overline{1,n},} $
$\alpha_i(t),$ $g(t) \in {\mathcal{C}}_{[0,T]}^1,$ functions $K_i(t,s)$
have continuous derivatives with respect to $t$ for $(t,s) \in cl({m_i}),$ $ K_n(t,t) \neq 0,$ $\alpha_i(0)=0,$ $0 < \alpha_1(t)<\alpha_2(t)<$ $\dots$ $< \alpha_{n-1}(t)<t,$ functions
$\alpha_1(t), \dots , \alpha_{n-1}(t) $  increase in a small neighborhood
$0\leq t \leq \tau,$
$0< \alpha_1^{\prime}(0) \leq$ $\dots$ $\leq \alpha_{n-1}^{\prime}(0)<1, $
$cl({m_i})$ denotes closure of set $m_i.$

Such weakly regular Volterra equations of the first kind with piecewise continuous kernels were first classified and generalised by Sidorov \cite{sid} and Lorenzi \cite{lor} and
extensively studied by many authors during the last decade. Here readers may refer to  monograph
\cite{bib4} and references therein. Volterra operator equations of the first kind were
studied by Sidorov and Sidorov \cite{sidsid}, sufficient conditions for existence of unique solution are obtained.  Tynda et al \cite{bib8} employed direct quadrature methods
for solution of  equations (\ref{math2021-e01}) both in linear and nonlinear cases. Muftahov and Sidorov \cite{muft} considered the numerical solution of nonlinear systems of such equations. Aghaei et al \cite{agh} applied
Legendre polynomials approximation method for solution of solution of linear Volterra integral equations with  piecewise continuous kernels. The numerical solution of the first kind Volterra convolution integral equations of the first kind  with
broad class of piecewise smooth kernels was considered by Davies and Duncan \cite{dav}. They employed the cubic convolution spline method and proved a stability bound.
Such Volterra models enjoys  applications in modeling various dynamical processes including storage systems \cite{eleng, matht}.
Generalized quadratures were employed by Sizikov and Sidorov \cite{siz} for solving singular
Volterra  integral equations of Abel type in application to infrared tomography.
The numerical solution of the
second-kind Volterra integral equation with weakly singular kernel is considered
in the piecewise polynomial collocation space by Linag and Brunner \cite{liang}. For conventional review of Volterra integral equations theory readers may refer e.g. to monographs by Brunner \cite{brun} and by Apartsyn \cite{apar}. Some studies of the Volterra integral equations of the first kind have led to the paradoxes as noted by Tynda in zbMATH~\cite{reft}.

In this paper, we also implement the CESTAC method and apply the CADNA library to find the numerical validation of the spline-collocation method to solve the problem (\ref{math2021-e01}). The priority of this strategy is to find the optimal step, accuracy and error of the numerical method. The paper is organised as following. In Section 2, the spline-collocation method is presented. Also the convergence of the method and smoothness of the solutions are studied. The use of Floating Point Arithmetic (FPA) is discussed in Section 3. The CESTAC method and its principle theorem can be found in Section 4. Using the this theorem we will show, how can we replace the conditions (\ref{cond1})   with (\ref{cond2}). The numerical results are illustrated in Section 5. Also in this section the comparative study between the results of the stochastic arithmetic (SA) and the FPA can be found.

%=====================================================================================================================
%=====================================================================================================================
%=====================================================================================================================
%=====================================================================================================================

\section{Polynomial spline-collocation method}
In this section, we construct a numerical method for solving problem \eqref{math2021-e01} -- \eqref{kernel}, based on the approximation of the exact solution by continuous local splines. First of all, within the framework of conditions \eqref{kernel}, we replace the original equation \eqref{math2021-e01} of the first kind with an equivalent equation of the second kind. To do this, we apply the standard technique of differentiating the equation:
\begin{equation}\label{math2021-e02}
   x(t)-\int_{0}^{t} h(t,s) x(s) ds = f(t), \; h(t,s)=-\frac{K'_{t}(t,s)}{K(t,t)}, \; f(t)=\frac{g'(t)}{K(t,t)}.
\end{equation}

Let us rewrite this equation in operator form
\begin{equation}\label{math2021-e03}
   (I-H)x=f, \; \text{where } (I-H)x\equiv x(t)-\int_{0}^{t} h(t,s) x(s) ds.
\end{equation}

%----------------------------------------
\subsection{Numerical scheme}

Let us introduce the partition of the interval $[0,T]$ with grid points $t_k, \; k=0,1,\ldots,N$. The introduced grid of nodes is not necessarily uniform and depends on the smoothness properties of the exact solution. Denote by $\Delta_k$ the segments $\Delta_k=[t_k,t_{k+1}],\; k=0,1,\ldots,N-1.$ Let then
\begin{equation}\label{Grid}
   \xi_k^j\in\Delta_k, \;  j=0,1,\ldots,r-1; \;\xi_k^0=t_k,\;\xi_k^{r-1}=t_{k+1}; \;k=0,1,\ldots,N-1,
\end{equation}
be additional nodes distributed in a certain way over the segment \(\Delta_k\). We denote by $P_r(x,\Delta_k)$ an operator, putting to the function $x (t),\;t\in\Delta_k,$  in correspondence an
interpolation polynomial of degree \(r-1\) for \(k=\overline{0,N-1}\) constructed on knots \(\xi_k^j.\)  Let then $x_N(t)$ be a local spline, defined on $[0,T]$ and composed of polynomials \(P_r(x,\Delta_k),\;k=0,1,\ldots,N-1.\)

We look for an approximate solution of \eqref{math2021-e03} as a spline \(x_N(t)\) with unknown values \(x_N(\xi_{k}^j)\),
\(k=0,1,\ldots,N-1, \; j=0,1,\ldots,r-1,\) at the knots of the grid.

The grid \eqref{Grid} depends on the considered class of functions to which the exact solutions belong and will be specified below.

The values \(x_N(\xi_{k}^i)\) in each segment \(\Delta_k, \; k=0,1,\ldots,N-1,\) are determined
step-by-step by the spline-collocation technique from the systems of linear equations
\begin{equation}\label{SplineCollocation}
  \begin{split}
   (I-H)P_N[x(t),\Delta_k]  \equiv P_N[x(t),\Delta_k]-\\
   -P_N\Biggl[\int\limits_{\Delta_k} P_N^{\tau}[h(t,s)]P_N[x(s),\Delta_k] ds,
    \;\Delta_k \Biggr]=P_N[f_k(t),\Delta_k].
  \end{split}
\end{equation}

Here $P_N$ is an operator of projection on the set of the local splines of the form \(x_N(t)\); \(f_k(t)\) is a new right part of
equation \eqref{math2021-e02} including the integrals over segments \(\Delta_j,\;j=0,1,\ldots,k-1,\) processed at the previous steps (in these domains, the spline parameters are already known):
\[
  f_k(t)=f(t)+\sum\limits_{l=0}^{k-1}\int\limits_{t_l}^{t_{l+1}}h(t,s)x_N(s)ds.
\]

%---------------------------------------------------------------------
\subsection{Convergence substantiation}
Let us rewrite equation \eqref{math2021-e02} and  projective method \eqref{SplineCollocation} in the operator form:
\begin{equation}\label{math2021-e04}
   x-Hx=f,\; H:X\to X,\; X\subset C(\Omega),\; \Omega=[0,T],
\end{equation}
\begin{equation}\label{math2021-e05}
  x_N-P_NHx_N=P_Nf,\;  P_N:X\to X_N,\; X_N\subset C(\Omega),
\end{equation}
where \(X\) is a dense set in \(C(\Omega)\) and \(X_N\) are the sets of corresponding local splines.

Since the homogenous Volterra integral equation \(x-Hx=0\) has only the trivial solution, the operator \(I-H\) is injective.
Hence, the operator \(I-H\) has the bounded inverse operator \((I-H)^{-1}:X\to X.\) For all sufficiently large \(N\) we have
the estimates
\[
  \|(I-P_NH)^{-1}\|_{C(\Omega)}=\|\Bigl((I-H)+(H-P_NH)\Bigr)^{-1}\|_{C(\Omega)} \leqslant
\]
\[
   \leqslant\frac{\|(I-H)^{-1}\|_{C(\Omega)}}{1-\|(I-H)^{-1}\|_{C(\Omega)}
   \|H-P_NH\|_{C(\Omega)}}
   \leqslant 2\|(I-H)^{-1}\|_{C(\Omega)}=A\;(const)
\]
if
\[
  \|H-P_NH\|_{C(\Omega)}\leqslant
  \frac{1}{2 \|(I-H)^{-1}\|_{C(\Omega)}}.
\]

Let us show that the last estimate holds for all sufficiently large \(N\). Since \(y(t)\equiv\bigl(Hx\bigr)(t)\in X\) and \(X\)
is a dense set in \(C(\Omega)\), we have
\[
  \|H-P_NH\|_{C(\Omega)}=
  \sup\limits_{x\in X,\|x\|\leqslant1}
  \max\limits_{t\in\Omega}|x(t)-P_Nx(t)|\leqslant e_N,
\]
where \( e_N\to 0\) as \(N\to\infty\). Therefore, \(\|H-P_NH\|_{C(\Omega)}\leqslant\frac{1}{2\|(I-H)^{-1}\|}\)
starting with sufficiently large \(N\).

Thus, the operators \((I-P_NH)^{-1}\) are exist and uniformly bounded and equation \eqref{math2021-e05} has a unique solution
for all sufficiently large \(N\) \cite{bib9}. Taking into account that \(P_Nx\to x\) as \(N\to\infty\) for all \(x\in X\), we apply the
projection operator \(P_N\) both to the left and the right parts of equation \eqref{math2021-e04}: $
  x-P_NHx=P_Nf+x-P_Nx.$  Subtracting this equation from \eqref{math2021-e05}, we obtain
$
  (I-P_NH)(x_N-x)=P_Nx-x,
$
and $(x_N-x)=(I-P_NH)^{-1}(P_Nx-x).$

This implies
\begin{equation}\label{math2021-Estimate}
    \|x_N-x\|_C\leqslant A\|P_Nx-x\|_C\leqslant  e_N(X).
\end{equation}

Thus, the accuracy of the approximate solution obtained via projective method \eqref{math2021-e05} is determined by the accuracy
\( e_N(X)\) of the approximation of functions from \(X\) by the local splines.

%=====================================================================================================================
%=====================================================================================================================
%=====================================================================================================================
%=====================================================================================================================

\subsection{Smooth solutions}

In this chapter, we describe in more detail the projection method \eqref{SplineCollocation} for the case of smooth input functions.
Namely, let the functions \(K_i(t,s), \; i=1,2,\ldots,n\) and \(g(t)\) in \eqref{math2021-e01}-\eqref{kernel} satisfy additional smoothness conditions.
Let \(g(t)\in C^{r+1}[0,T]\), \(K_i(t,s)\in C^{r+1,r}[0,T]^2\) (continuously differentiable on each variable).
The exact solution \(x(t)\) of the equation \eqref{math2021-e01} in this case belongs to \(X=C^r[0,T]\)  \cite{bib4}.

We introduce the uniform partition of the interval $[0,T]$ with grid points $t_k=\frac{kT}{N}, \; k=0,1,\ldots,N$.
Denote by $\Delta_k$ the segments $\Delta_k=[t_k,t_{k+1}],\; k=0,1,\ldots,N-1.$ Let
$
  \xi_k^j=\frac{t_{k+1}+t_k}{2}+\frac{t_{k+1}-t_k}{2}y_j,
$ and
$
  j=1,2,\ldots,r-2;\;\xi_0^k=t_k,\;\xi_{r-1}^k=t_{k+1};
  \;k=0,2,\ldots,N-1.
$
where $y_j$ are the roots of the Legendre polynomials of degree \(r-2\).

We denote by $P_r(x,\Delta_k)$ an operator, putting to the function $x (t),\;t\in\Delta_k $  in correspondence an
interpolation polynomial of degree \(r-1\) for \(k=\overline{0,N-1}\) constructed on knots \(\xi_k^j.\) Let then $x_N(t)$ be a local spline, defined on $[0,T]$ and composed of polynomials \(P_r(x,\Delta_k),\;k=0,1,\ldots,N-1.\) We look for an approximate solution of the equation \eqref{math2021-e01} as a spline $x_N(t),$ $0\leqslant t\leqslant T,$
with unknown coefficients \(x_k^j,\) \(k=\overline{0,N-1}.\) Let us describe the process of definition \(x_k^j\). At first we find the coefficients \(x_0^j,\;j=0,1,\ldots,r-1,\) from the system of equations
\begin{equation}\label{math2021-e06}
   P_r(x,\Delta_0)(\xi_0^j)-P_r\left[\int\limits_{0}^{\xi_0^j}
   h(\xi_0^j,\tau)P_r(x,\Delta_0)(\tau)d\tau,\Delta_0
   \right]=P_r(f,\Delta_0)(\xi_0^j).
\end{equation}
Here for integrals calculation in \eqref{math2021-e06} we employ compound Gaussian quadrature rule with \(r-2\) points constructed on the auxiliary mesh
linked to the lines  \(\alpha_i(t),\; i=\overline{1,n},\) of the kernel discontinuities for each specific value of \(N\).
Note also that all values of the unknown function in intermediate points are computed with help of interpolation.
The resulting system of linear algebraic equations is then solved by the Jordan-Gauss method.

In order to determine the coefficients $x_1^j,\;j=0,1,\ldots,r-1,$ of the local spline $x_N(t)$ on the segment $\Delta_1$ we
represent \eqref{math2021-e02} in the following form
\begin{equation}\label{math2021-e07}
    x(t)-\int\limits_{\xi_1^0}^{t}h(t,\tau)x(\tau)d\tau=f_1(t), \; f_1(t)=f(t)+\int\limits_{0}^{\xi_0^{r-1}}h(t,\tau)P_r(x(\tau),\Delta_0)d\tau.
\end{equation}
The equation \eqref{math2021-e07} is then solved by analogy, using the scheme \eqref{math2021-e06}.  Repeating this process  $N$ times we obtain the approximate solution \(x_N(t)\) of the equation \eqref{math2021-e02} over all interval \([0,T].\)

%**********************
The error of approximation of the exact solution by the polynomials constructed in this way at each step of the process can be estimated by the following inequality (here readers may refer to \cite{bib7})
\begin{equation}\label{LocalError}
  \|x(t)-x_N(t)\|_{C[\Delta_k]}\leq \frac{L_{r}\left(\frac{T}{N}\right)^{r}}{r!}, \; k=0,1,\ldots,N-1,
  \quad L_r=\max\limits_{t\in\Delta_k}\left|x^{r}(t)\right|.
\end{equation}

Taking into account the general estimate \eqref{math2021-Estimate} of the error of the method, we obtain the error estimate in this case (\(X=C^r[0,T]\))
\begin{equation}\label{SplineError}
    \|x(t)-x_N(t)\|_{C[0,T]}\asymp N^{-r}.
\end{equation}
%**********************
Boikov and Tynda  \cite{bib10,bib11}  established that such type numerical methods for Volterra integral equations are also optimal with respect to complexity and accuracy order. Thus, an effective projective method for solving equations of the form \eqref{math2021-e01} is proposed.
%=====================================================================================================================
%=====================================================================================================================
%=====================================================================================================================
%=====================================================================================================================
\section{Using the Floating Point Arithmetic}

In general form when  floating point arithmetic (FPA) is employed  it is necessary to have the exact and approximate solutions $x(t)$ and $x_N(t)$ and also small value $\epsilon$ to use the following conditions
\begin{equation}\label{cond1}
  |x(t) - x_N(t)| \leq \epsilon, ~~~or ~~~  |x_N(t) - x_{N-1}(t)| \leq \epsilon.
\end{equation}
But the main problem is that the exact solution and optimal $\epsilon$ are unknowns. Thus by putting small values instead of $\epsilon$, the approximate solution will not be accurate and for large values we will have many extra iterations. In order to avoid these problems, the CESTAC method and the CADNA library will be utilized \cite{s29}. In this novel method, the accuracy will be obtained using successive iterations $x_N(t)$ and $x_{N-1}(t)$ and the following condition
\begin{equation}\label{cond2}
\begin{array}{l}
  |x_N(t) - x_{N-1}(t)| = @.0.
\end{array}
\end{equation}
We apply this condition to control the accuracy of the method and avoid extra iterations by using number of common significant digits (NCSDs) between $x_N(t)$ and $x_{N-1}(t)$. $@.0$ in Eq. (\ref{cond2}) displays the informatical zero which can be produced only in the CESTAC method by the CADNA library. It shows  that the NCSDs of two successive approximations and approximate and exact solutions are almost equal to zero. Vignes and La Porte \cite{s32} have presented the method for the first time in 1974. In \cite{s29} Vignes has described the CESTAC method to evaluate the numerical results of some computational methods. Some conditions of the CESTAC method, applying different tools to write the CADNA codes \cite{s26} and also some properties of the SA have been studied by  Chesneaux. All the CESTAC evaluations should be accomplished using the CADNA library. Handling this scheme, the optimal results, step and error of the method can be recognized.   Lamotte et al. has implemented the CESTAC method using C and C++ codes.  J\'{e}z\'{e}quel et al have discussed the new version of the CADNA library using Fortran programs. Recently applying this method to control the accuracy of the Taylor expansion method to solve the generalized Abel's integral equation \cite{s2}, mathematical model of Malaria infection \cite{s4,s6}, nonlinear fractional order model of COVID-19 \cite{s5}, solving nonlinear shallow water wave equation \cite{s7}, Adomian decomposition method, homotopy perturbation method and Taylor-collocation method  for solving Volterra integral equation \cite{s8,s11,s15, taylorisu}, dynamical control of the reverse osmosis system \cite{s9,s10}, solving integrals using the numerical methods have been done. Moreover the CESTAC method has been used to find the optimal convergence control parameter of the homotopy analysis method in both fuzzy and crisp forms \cite{s18,s19}.

\section{CESTAC Method}

The CESTAC method is a powerful and applicable tool to validate the numerical results of numerical procedures. It should be applied based on the SA. Let $B$ be a set of reproduced values by  computer. For real value $g^{*}$, we can find a member of set $B$ such as $G^{*}$ with $\alpha$ mantissa bits of the binary FPA as $
G^{*}=g^{*}-\rho 2^{E-\alpha}\phi,
$
where the sign showed by $\rho$, the missing segment of the mantissa presented by $2^{-\alpha}\phi$ and the binary exponent of the
result displayed by $E$. Replaying 24 and 53 instead of $\alpha$, the results can be found by single and double accuracies. Assuming $\phi$ as a stochastic variable and having uniformly distribution on $[-1, 1]$, we will be able to make perturbation on the last mantissa bit of $g^{*}$. Thus for the obtained results of $G^{*}$, the mean and standard deviation values  $(\mu)$ and $(\sigma)$ can be found. Doing the mentioned scheme $p$-times $p$ samples of $G^{*}$ can be produced as
$ \Phi = \left\{ G^{*}_1, G^{*}_2, ..., G^{*}_p\right\}.$
Thus the mean and standard deviation can be found as follows
$\displaystyle \tilde{G}^{*} = \frac{\sum_{k=1}^p
G^{*}_k}{p},~~~~~\displaystyle \sigma^2= \frac{\sum_{k=1}^p
( G^{*}_k - \tilde{G}^{*} )^2}{p-1}.$

Using the mentioned computations the NCSDs of $G^{*}$ and $\tilde{G}^{*}$ can be generated using the following relation
$\displaystyle \mathcal{C}_{\tilde{G}^{*}, G^{*}} = \log_{10}\frac{\sqrt{p}
\left|\tilde{G}^{*} \right|}{\tau_{\delta} \sigma},$
where $\tau_{\delta}$ is the value of $T$ distribution as the confidence
interval is $1-\delta$, with $p-1$  freedom degree. Showing $G^{*}=@.0$, the process stopped if we have $\tilde{G}^{*} =0,$ or $\mathcal{C}_{\tilde{G}^{*}, G^{*}}
\leq 0.$

In this method, the mathematical softwares Mathematica, Maple or MATLAB must be replaced by the  CADNA library. This library should be implemented on the LINUX operating system and we all codes should be compiled by  C, C++, FORTRAN or ADA. The main benefit of the method is to find the optimal results, step size and error of the method.

%=============================================================================================================================

\begin{definition} \label{def} For $\ell_1, \ell_2 \in \mathbb{R}$  the NCSDs can be defined as
\begin{equation}\label{8}
\mathcal{C}_{\ell_1,\ell_2}=\left\{
\begin{array}{l}
\  \log_{10}\left|\frac{\ell_1+\ell_2}{2(\ell_1-\ell_2)} \right| =
\log_{10}\left|\frac{\ell_1}{\ell_1-\ell_2} - \frac{1}{2} \right|,~~~~~\ell_1\neq \ell_2, \\
    \\
    \infty,~~~~~otherwise.
\end{array}
\right.
\end{equation}
\end{definition}

\begin{theorem} \label{th4}
Assume that $x(t)$ and $x_N(t)$ are the exact and approximate solutions of Eq. (\ref{math2021-e01}). The NCSDs of two successive approximations are almost equal to the NCSDs of exact and approximate solutions and we have
\begin{equation}\label{9}
\displaystyle   \mathcal{C}_{x_{N}(t),x_{N+1}(t)} \simeq  \mathcal{C}_{x_{N}(t),x(t)}.
\end{equation}
\end{theorem}

\textbf{Proof:} Applying Definition 1 and Eq. (\ref{SplineError}) for to iterations $x_N(t)$ and $x_{N+1}(t)$ we can write
$$
\begin{array}{l}
\displaystyle \mathcal{C}_{x_N(t),x_{N+1}(t)} =\log_{10}\left|
\frac{x_N(t)+x_{N+1}(t)}{2(x_N(t)-x_{N+1}(t))}\right| =
\log_{10}\left|
\frac{x_N(t)}{x_N(t)-x_{N+1}(t)} - \frac{1}{2} \right| \\
\\
~~~~~~~~~~~~~~~~~~~~\displaystyle
=\log_{10}\left|\frac{x_N(t)}{x_N(t)-x_{N+1}(t)} \right| +
\log_{10}\left|1- \frac{1}{2
x_N(t)}(x_N(t)-x_{N+1}(t)) \right|\\
\\
~~~~~~~~~~~~~~~~~~~~\displaystyle
=\log_{10}\left|\frac{x_N(t)}{x_N(t)-x_{N+1}(t)} \right| +
\mathcal{O} \left(x_N(t)-x_{N+1}(t)\right).
\end{array}
$$

We know
$
x_N(t)-x_{N+1}(t) = x_N(t)-x(t)-(x_{N+1}(t)-x(t)) = E_n (t) -
E_{n+1}(t),
$
therefore we get
$
 \mathcal{O}\left(x_N(t)-x_{N+1}(t)\right) =  \mathcal{O}\left( N^{-r} \right) +  \mathcal{O}\left( (N+1)^{-r} \right)= \mathcal{O}\left( N^{-r} \right).
$
And finally we have
\begin{equation}\label{15}
\mathcal{C}_{x_N(t),x_{N+1}(t)} =
\log_{10}\left|\frac{x_N(t)}{x_N(t)-x_{N+1}(t)} \right|
+\mathcal{O}\left(N^{-r}\right).
\end{equation}
Repeating the process for exact and approximate solutions, the following relation can be obtained
\begin{equation}\label{16}
 \begin{array}{ll}
\displaystyle \mathcal{C}_{x_N(t),x(t)} &\displaystyle = \log_{10}\left|
\frac{x_N(t)+x(t)}{2(x_N(t)-x(t))}\right| = \log_{10}\left|
\frac{x_N(t)}{x_N(t)-x(t)} - \frac{1}{2} \right| \\
\\
&\displaystyle =\log_{10}\left| \frac{x_N(t)}{x_N(t)-x(t)}\right| +
\mathcal{O}(x_N(t) -
x(t))\\
\\
&\displaystyle = \log_{10}\left| \frac{x_N(t)}{x_N(t)-x(t)}\right|
+\mathcal{O}\left(N^{-r}\right).
\end{array}
\end{equation}

Based on  Eqs. (\ref{15}) and (\ref{16}) we can write
$$\begin{array}{ll}
\displaystyle \mathcal{C}_{x_N(t),x(t)} - \mathcal{C}_{x_N(t),x_{N+1}(t)} &\displaystyle
= \log_{10}\left| \frac{x_N(t)}{x_N(t)-x(t)}\right| -
\log_{10}\left|\frac{x_N(t)}{x_N(t)-x_{N+1}(t)}
\right| \\
\\
&\displaystyle +\mathcal{O}\left(N^{-r}\right) = \log_{10}\left| \frac{x_N(t)-x(t)
}{x_N(t)-x_{N+1}(t) } \right| +\mathcal{O}\left(N^{-r}\right)\\
\\
&\displaystyle = \log_{10}\left|
\frac{\mathcal{O}\left(N^{-r}\right)}{\mathcal{O}\left(N^{-r}\right)} \right| +\mathcal{O}\left(N^{-r}\right) = \mathcal{O}\left(N^{-r}\right)\\
  \end{array}
$$
and
$
\displaystyle \mathcal{C}_{x_N(t),x(t)} - \mathcal{C}_{x_N(t),x_{N+1}(t)} =
\mathcal{O}\left(N^{-r}\right).
$
Clearly by approaching $N$ to infinity, $\mathcal{O}\left(N^{-r}\right)$ tends zero and we obtain
$
\mathcal{C}_{x_N(t),x(t)} = \mathcal{C}_{x_N(t),x_{N+1}(t)}.
$

%
% ==========Numerical results ==========================
%=====================================================================================================================
%=====================================================================================================================
%=====================================================================================================================
%=====================================================================================================================
\section{Numerical results}

To illustrate the effectiveness of the suggested polynomial spline-collocation method, we present the results for two test equations.

\textbf{Example 1.}
Consider the following integral equation
\begin{equation}\label{Model-01}
  \int\limits_0^{\frac{t}{2}}(t+s)x(s)ds+\int\limits_{\frac{t}{2}}^{\frac{2t}{3}} t s x(s)ds+\int\limits_{\frac{2t}{3}}^{t}e^sx(s)ds=f(t), \quad t \in [0; 1],
\end{equation}
where the right side of the equation was chosen so that the exact solution was \(x^*=t\sin t\).
The following designations are used in the tables below: \(\mathbf{N}\) is the number of segments of the main partition, \(\mathbf{r}\) is the parameter responsible for the order of the spline,  \( e=||x_N(t)-x^*(t)||_{C_{[0, T]}}\).
\begin{center}
  \begin{table}[h]
   \caption{The error for \eqref{Model-01} at the value \(\mathbf{r}=4\).  }\label{Tt01}
    \centering\scalebox{0.7}{
  \begin{tabular}{|c|c|c|c|c|c|c|c|} \hline
 $\mathbf{N}$ & 1 & 5 & 10 & 20 & 50 & 100 & 500   \\  \hline
 $ e$ & 6.57\(\cdot 10^{-4}\)  &  2.38\(\cdot 10^{-7}\)  & 7.55\(\cdot 10^{-9}\) & 2.38\(\cdot 10^{-10}\) & 2.45\(\cdot 10^{-12}\) & 7.65\(\cdot 10^{-14}\) & 2.46\(\cdot 10^{-17}\)  \\  \hline
  \end{tabular}}
  \end{table}
\end{center}
%*************************************
\begin{center}
  \begin{table}[h]
    \caption{The error for \eqref{Model-01} at the value \(\mathbf{r}=7\).  }\label{Tt02}
    \centering\scalebox{0.7}{
  \begin{tabular}{|c|c|c|c|c|c|c|c|} \hline
 $\mathbf{N}$ & 1 & 5 & 10 & 20 & 50 & 100 & 500   \\  \hline
 $ e$ & 2.95\(\cdot 10^{-7}\)  &  4.71\(\cdot 10^{-12}\)  & 3.73\(\cdot 10^{-14}\) & 2.94\(\cdot 10^{-16}\) & 4.82\(\cdot 10^{-19}\) & 3.77\(\cdot 10^{-21}\) & 1.39\(\cdot 10^{-25}\)  \\  \hline
  \end{tabular}}
  \end{table}
\end{center}
%*****************************************
\begin{figure}[h]
    \center{\includegraphics[width=0.75\linewidth]{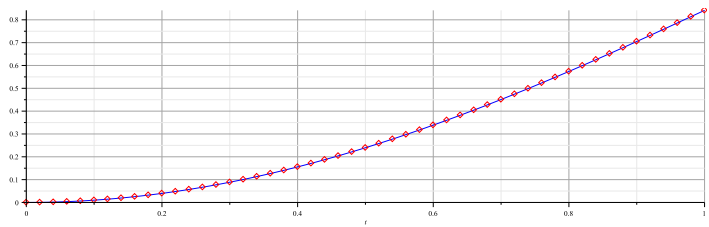}}
    \caption{The exact and approximate solution of \eqref{Model-01} with \(N=5,\; r=8\).}
    \label{Figure1}
\end{figure}
%*****************************************

Clearly, Tables \ref{Tt01} and \ref{Tt02} depend on  the exact solution. For Table \ref{tb1} we apply the CESTAC method and the results are obtained $\textbf{r}=5$ and the optimal results are
$
  N_{opt}=6,
  x_6(0.05)= 0.2498959,
  error_{opt} = 0.2328306E-009.
$
Table \ref{tbf1} is obtained using the spline-collocation method and the FPA for the same value of  $\textbf{r}$. It is obvious that for  $\epsilon=10^{-2}$ the algorithm will be stopped at $N=1$ and for $\epsilon = 10^{-10}$ we have $N=4$. For small values of $\epsilon$ we will need to provide many extra iterations without improving the accuracy. Tables \ref{tb2} and \ref{tbf2} present the results for $\textbf{r}=10$ using the CESTAC method and we have
$  N_{opt}=2,
  x_6(0.05)= 0.2498958.
$
It means that we do not need to produce a smaller partition and we can stop at the specified value of $N$. Thus according to the obtained results, $N=2$ is enough and we do not need to find more results.

 \begin{table}[h!]
\caption{ Results of the CESTAC method for $\textbf{r}=5$.  }\label{tb1} \centering\scalebox{0.7}{
\begin{tabular}{|l|l|l|l|}
  \hline
$N$ ~~~~~~~~~& $x_{N+1}(t)$~~~~~~~~~ &$|x_{N+1}(t) - x_{N}(t)|$ ~~~~~~~~~&$|x_{N}(t) - x^*(t)|$  ~~~~~~~~~\\
    \hline
    1     &       0.2474730     &        0.2474730E-002  &      0.24227E-004\\
    2     &       0.2498547     &        0.23816E-004   &     0.411E-006\\
    3     &       0.2498938     &        0.391E-006     &   0.19E-007\\
    4     &       0.2498958     &        0.19E-007     &   0.2328306E-009\\
    5     &       0.2498959     &        0.6984919E-009   &     @.0\\
    6     &       0.2498959     &        @.0               &            @.0\\
\hline
  \end{tabular}}
\end{table}

%Example 1 (for r=5):

%N	x(t) - exact	x_N(t) -approx	Error at the point t=0.05

\begin{table}[h!]
\caption{ The spline-collocation using the FPA for solving Example 1 with $\textbf{r}=5$.  }\label{tbf1} \centering\scalebox{0.7}{
\begin{tabular}{|l|l|l|}
  \hline
$N$ & $x_{N+1}(t)$  &$|x_{N+1}(t) - x^*(t)|$  \\
    \hline
1	&	0.0024747309902353909851940001084329&	0.0000242274732985254545492499338416\\
2	&	0.0024985471725248365458815553649446&	4.112910090798938616946773299e-7\\
3	&	0.0024989386743724730905507916941931&	1.97891614433491924583480814e-8\\
4	&	0.0024989586478083542254432950394024&	1.842744377857000449971279e-10\\
5	&	0.0024989592632306142870868948525715&	7.996966978473436448102970e-10\\
6	&	0.0024989588352407427373351401205658&	3.717068262975918900782913e-10\\
7	&	0.0024989585878160274749598471153271&	1.242821110352165970730526e-10\\
8	&	0.0024989584860466374183598209171807&	2.25127209786165708749062e-11\\
9	&	0.0024989584529530046839100658673684&	1.05809117558331841749061e-11\\
10	&	0.0024989584476966830750021803651006&	1.58372333647410696771739e-11\\
\hline
  \end{tabular}}
\end{table}

 \begin{table}[h!]
\caption{ Results of the CESTAC method for $\textbf{r}=10$.  }\label{tb2} \centering\scalebox{0.7}{
\begin{tabular}{|l|l|l|l|}
  \hline
$N$ ~~~~~~~~~& $x_{N+1}(t)$ ~~~~~~~~~&$|x_{N+1}(t) - x_{N}(t)|$ ~~~~~~~~~&$|x_{N+1}(t) - x^*(t)|$ ~~~~~~~~~ \\
    \hline
        1    &        0.2498958    &         0.2498958E-002  &      @.0\\
    2        &    0.2498958       &      @.0                 &          @.0\\
\hline
  \end{tabular}}
\end{table}

%Example 1 (for r=10):

%N	x(t) - exact	x_N(t) -approx	Error at the point t=0.05

\begin{table}[h!]
\caption{ The spline-collocation using the FPA for solving Example 1 with $\textbf{r}=10$. }\label{tbf2} \centering\scalebox{0.7}{
\begin{tabular}{|l|l|l|}
  \hline
$N$ & $x_{N+1}(t)$  &$|x_{N+1}(t) - x^*(t)|$  \\
    \hline
1&		0.0024989584688525888359310360465249&	5.3186723961877860042504e-12\\
2&		0.0024989584635324124519264394003959&	1.5039878168106418786e-15\\
3&		0.0024989584635338888159300590349018&	2.76238131910073727e-17\\
4	&	0.0024989584635339165728796529239866&	1.331364028817121e-19\\
5&		0.0024989584635339164599105607342320&	2.01673106919575e-20\\
6&		0.0024989584635339164307566189869647&	8.9866310553098e-21\\
7&		0.0024989584635339164379496352522859&	1.7936147899886e-21\\
8&		0.0024989584635339164399427718497446&	1.995218074701e-22\\
9	&	0.0024989584635339164398662187644619&	1.229687221874e-22\\
10	&	0.0024989584635339164397061697372937&	3.70803049808e-23\\
\hline
  \end{tabular}}
\end{table}

%=====================================================================================================================
%=====================================================================================================================
%=====================================================================================================================
%=====================================================================================================================

\textbf{Example 2.}
Consider the following integral equation
\begin{equation}\label{Model-02}
  \int\limits_0^{\frac{t}{3}}(t-s)^2x(s)ds+\int\limits_{\frac{t}{3}}^{\frac{3t}{4}} \cos(s) x(s)ds+\int\limits_{\frac{3t}{4}}^{t}(1+\sin(2s))x(s)ds=f(t), \quad t \in [0; 2],
\end{equation}
where the right side of the equation was chosen so that the exact solution was \(x^*=e^{2-t}t^2\).
\begin{center}
  \begin{table}[h]
   \caption{The error for \eqref{Model-02} at the value \(\mathbf{r}=5\).  }\label{math2021-t03}
    \centering\scalebox{0.7}{
  \begin{tabular}{|c|c|c|c|c|c|c|c|} \hline
 $\mathbf{N}$ & 1 & 5 & 10 & 20 & 50 & 100 & 500   \\  \hline
 $ e$ & 7.67\(\cdot 10^{-3}\)  &  4.89\(\cdot 10^{-6}\)  & 1.70\(\cdot 10^{-7}\) & 5.61\(\cdot 10^{-9}\) & 5.96\(\cdot 10^{-11}\) & 1.88\(\cdot 10^{-12}\) & 6.09\(\cdot 10^{-16}\)  \\  \hline
  \end{tabular}
}
  \end{table}
\end{center}
%*************************************
\begin{center}
  \begin{table}[h]
   \caption{The error for \eqref{Model-02} at the value \(\mathbf{r}=10\).  }\label{math2021-t04}
    \centering\scalebox{0.7}{
  \begin{tabular}{|c|c|c|c|c|c|c|c|} \hline
 $\mathbf{N}$ & 1 & 5 & 10 & 20 & 50 & 100 & 500   \\  \hline
 $ e$ & 8.61\(\cdot 10^{-9}\)  &  1.41\(\cdot 10^{-15}\)  & 1.46\(\cdot 10^{-18}\) & 1.47\(\cdot 10^{-21}\) & 1.49\(\cdot 10^{-23}\) & 1.43\(\cdot 10^{-24}\) & 7.47\(\cdot 10^{-27}\)  \\  \hline
  \end{tabular}
}
  \end{table}
\end{center}
%*****************************************
\begin{figure}[h]
    \center{\includegraphics[width=0.75\linewidth]{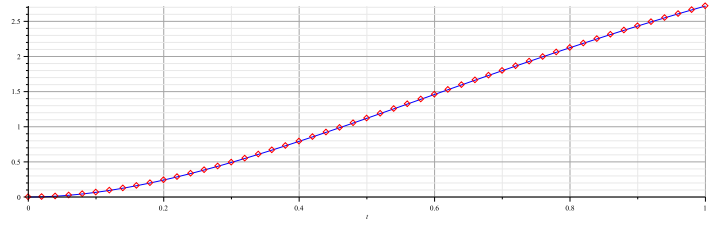}}
    \caption{The exact and approximate solution of \eqref{Model-02} with \(N=5,\; r=8\).}
    \label{Figure2}
\end{figure}
%*****************************************

All calculations were performed in the Maple system with parameter \texttt{Digits:=30;} (the number of digits that Maple uses when making calculations with software floating-point numbers). As we can see from Tables \ref{math2021-t03} and \ref{math2021-t04}, the practical error of the method corresponds to the theoretical estimate \eqref{SplineError}. All the results of Tables \ref{tb3} and \ref{tb4}, are obtained using the CESTAC method. For $\textbf{r}=6$ we get
$
  N_{opt}=6,
  x_6(0.05)= 0.1757171,
  error_{opt} = 0.9313225E-008,
$
and for $\textbf{r}=12$ we have
$  N_{opt}=2,
  x_6(0.05)= 0.1757171.
$
According to the results for large values of  $\textbf{r}$ the results are more accurate. Tables \ref{tbf3} and \ref{tbf4} are obtained for the spline-collocation method using the FPA. Comparing the results of  the FPA and the SA, we can introduce the CESTAC method as a good tool to control the accuracy and the step size of the spline-collocation method for solving the mentioned problem.

 \begin{table}[h!]
\caption{Results of the CESTAC method for $\textbf{r}=6$. }\label{tb3}
\centering\scalebox{0.7}{
\begin{tabular}{|l|l|l|l|}
  \hline
$N$ ~~~~~~~~~& $x_{N+1}(t)$~~~~~~~~~ &$|x_{N+1}(t) - x_{N}(t)|$ ~~~~~~~~~&$|x_{N}(t) - x^*(t)|$ ~~~~~~~~~ \\
    \hline
    1      &      0.1762535     &        0.1762535E-001&        0.5363E-004\\
    2      &      0.1757219    &         0.5316E-004   &     0.47E-006\\
    3      &      0.1757169     &        0.49E-006    &    0.2E-007\\
    4      &      0.1757170     &        0.1E-007     &   0.9313225E-008\\
    5      &      0.1757171     &        0.7450580E-008  &      @.0\\
    6      &      0.1757171     &        @.0             &              @.0\\
\hline
  \end{tabular}}
\end{table}

 \begin{table}[h!]
\caption{ Results of the CESTAC method for $\textbf{r}=12$. }\label{tb4} \centering\scalebox{0.7}{
\begin{tabular}{|l|l|l|l|}
  \hline
$N$ ~~~~~~~~~& $x_{N+1}(t)$ ~~~~~~~~~&$|x_{N+1}(t) - x_{N}(t)|$ ~~~~~~~~~&$|x_{N}(t) - x^*(t)|$ ~~~~~~~~~ \\
    \hline
    1      &      0.1757171    &         0.1757171E-001  &      @.0\\
    2      &      0.1757171    &         @.0              &             @.0\\
\hline
  \end{tabular}}
\end{table}

%Example 2 (for r=6):

%N	x(t) - exact	x_N(t) -approx	Error at the point t=0.05

\begin{table}[h!]
\caption{  The spline-collocation using the FPA for solving Example 2 with $\textbf{r}=6$. }\label{tbf3} \centering\scalebox{0.7}{
\begin{tabular}{|l|l|l|}
  \hline
$N$ & $x_{N+1}(t)$  &$|x_{N+1}(t) - x^*(t)|$  \\
    \hline
1	& 	0.017625357759763366370583965414699&	0.000053638808290133034856760580788\\
2	& 	0.017572191765083977687409203731931&	4.72813610744351681998898020e-7\\
3	& 	0.017571692863117541262788927971815&	2.6088355692072938276862096e-8\\
4	& 	0.017571708598170673074341558874276&	1.0353302560261385645959635e-8\\
5	& 	0.017571717076595299740419367871963&	1.874877933595307836961948e-9\\
6	& 	0.017571719104941762543174921162397&	1.53468529207447716328486e-10\\
7	& 	0.017571719329605570787973205575224&	3.78132337452246000741313e-10\\
8	& 	0.017571719187368398597882888516563&	2.35895165262155683682652e-10\\
9	& 	0.017571719040773236192761953789278&	8.9300002857034748955367e-11\\
10	& 	0.017571718954164366107614229639737&	2.691132771887024805826e-12\\
\hline
  \end{tabular}}
\end{table}

%Example 2 (for r=12):

%N	x(t) - exact	x_N(t) -approx	Error at the point t=0.05

\begin{table}[h!]
\caption{ The spline-collocation using the FPA for solving Example 2 with $\textbf{r}=12$. }\label{tbf4} \centering\scalebox{0.7}{
\begin{tabular}{|l|l|l|}
  \hline
$N$ & $x_{N+1}(t)$  &$|x_{N+1}(t) - x^*(t)|$  \\
    \hline
1	&	0.017571718951506900943716697361262&	3.3667607989492527351e-14\\
2	&	0.017571718951473218595173777454132&	1.4740553427379779e-17\\
3	&	0.017571718951473233466169941902544	&1.30442737068633e-19\\
4	&	0.017571718951473233336814972540526	&1.087767706615e-21\\
5	&	0.017571718951473233335480404651884&	2.46800182027e-22\\
6	& 	0.017571718951473233335711266248458&	1.5938585453e-23\\
7	& 	0.017571718951473233335728475858345&	1.271024434e-24\\
8	& 	0.017571718951473233335727478289445&	2.73455534e-25\\
9	& 	0.017571718951473233335727218074944&	1.3241033e-26\\
10	& 	0.017571718951473233335727227475141	&2.2641230e-26\\
\hline
  \end{tabular}}
\end{table}

\subsection{Stability experiments}

To illustrate the stability of suggested numerical method, we introduced a random error in calculating the values of free term \(f(t)\) of the equations \eqref{Model-01} and \eqref{Model-02}.
The range of introduced random errors is \((-\delta,\delta)\).

Depending on the \(\delta\), the following results are obtained.
\begin{center}
  \begin{table}[h!]\label{Perturbed1}
    \centering\scalebox{0.7}{
  \begin{tabular}{|c|c|c|c|c|c|c|c|} \hline
  ~~~$\delta$  ~~~& ~~~0 ~~~& ~~~$10^{-6}$ ~~~& ~~~$10^{-5}$ ~~~& ~~~$10^{-4}$ ~~~& ~~~$10^{-3}$ ~~~& ~~~$10^{-2}$ ~~~  \\  \hline
 $\varepsilon$ & 7.17\(\cdot 10^{-7}\)  &  2.74\(\cdot 10^{-5}\)  & 2.75\(\cdot 10^{-4}\) & 0.00275 & 0.0275 & 0.2746  \\  \hline
  \end{tabular}}
  \caption{The error for \eqref{Model-01} with \(\mathbf{N}=5,\; \mathbf{r}=5\).  }
  \end{table}
\end{center}
%*************************************
\begin{center}
  \begin{table}[h!]\label{Perturbed2}
    \centering\scalebox{0.7}{
  \begin{tabular}{|c|c|c|c|c|c|c|c|} \hline
 ~~~ $\delta$ ~~~& ~~~0~~~ & ~~~$10^{-6}$ ~~~& ~~~$10^{-5}$ ~~~& ~~~$10^{-4}$ ~~~& ~~~$10^{-3}$ ~~~& ~~~$10^{-2}$  ~~~ \\  \hline
 $\varepsilon$ & 4.89\(\cdot 10^{-6}\)  &  5.71\(\cdot 10^{-5}\)  & 5.28\(\cdot 10^{-4}\) & 0.00524 & 0.05235 & 0.52351  \\  \hline
  \end{tabular}}
  \caption{The error for \eqref{Model-02} with \(\mathbf{N}=5,\; \mathbf{r}=5\).  }
  \end{table}
\end{center}
%*************************************
From the results proposed in the tables 13 and 14 %\ref{Perturbed1} and \ref{Perturbed2},
it is possible to judge the continuous dependence of the solution on the initial data and conclude about the stability of the numerical method.
This result is not surprising: despite the fact that the initial equation is an equation of the first kind, it is solved in such spaces in which the problem is well-posed.
The case of an ill-posedness (when noisy initial data may lead to instability and regularization will be required) will be studied in future works.

\bigskip

\section{Conclusion}

We have applied the spline-collocation method for solving the Volterra integral equations of the first kind with discontinuous kernel. The convergence of the method and the smoothness of the solution have been discussed. Using the CESTAC method we have tried to control the accuracy and step size of the method. The principle theorem of the CESTAC method will help us to apply the condition (\ref{cond2}) instead of (\ref{cond1}). Thus we will be able to find the optimal results, optimal error and optimal step of the method.

\bigskip
%пїЅпїЅпїЅпїЅпїЅпїЅпїЅпїЅ пїЅпїЅ пїЅпїЅпїЅпїЅпїЅпїЅ пїЅпїЅ пїЅпїЅпїЅпїЅпїЅ пїЅпїЅпїЅпїЅпїЅпїЅ

\textbf{Aleksandr Tynda}, Candidate of Science  (Physics and Mathematics), Professor,
Higher and Applied Mathematics Department, Penza State University, 440026 Penza, Russia;\\
email: tyndaan@mail.ru, \\ORCID ID https://orcid.org/0000-0001-6023-9847.

\textbf{Samad Noeiaghdam}, PhD (Applied Mathematics), Associate Professor,
Industrial Mathematics Laboratory, Baikal School of BRICS, Irkutsk National Research Technical University, Irkutsk, 664074, Russia; \\
Department of Applied Mathematics and Programming, South Ural State University, Lenin prospect 76, Chelyabinsk, 454080, Russia; \\
 email: snoei@istu.edu; noiagdams@susu.ru;\\  ORCID ID https://orcid.org/0000-0002-2307-0891.

\textbf{Denis Sidorov}, Doctor of Science (Physics and Mathematics), Professor of RAS,
Energy Systems Institute, Siberian Branch of Russian Academy of Science, 664033 Irkutsk, Russia\\
Industrial Mathematics Laboratory, Baikal School of BRICS, Irkutsk National Research Technical University, Irkutsk, 664074, Russia; \\
Institute of Mathematics and Information Technologies, Irkutsk State University, Irkutsk 664025, Russia;
 email: sidorovdn@istu.edu,\\ ORCID ID https://orcid.org/0000-0002-3131-1325.

%
%\end{article}
%
%\end{document}

%название
\naze{Метод полиномиальной сплайн-коллокации для решения слабо регулярных интегральных уравнений Вольтерра I рода}

%авторы

%\avtore{А. Тында$^{1}$, П. П. Петров$^{2}$}
%
%\inst{\small {$^1$ Иркутский государственный университет, Иркутск, Российская Федерация\\
%$^2$ Новосибирский государственный университет, Новосибирск, Российская Федерация}}

\avtore{A. Тында $^1$,  С. Нойягдам$^{2,3}$, Д. Сидоров$^{4,2,5}$}

\noindent \inst{\small {$^1$ Кафедра высшей и прикладной математики, Пензенский государственный университет,  Пенза, Россия\\
$^2$ Лаборатория промышленной математики, Байкальская школа БРИКС,
Иркутский национальный исследовательский университет, Иркутск,  Россия\\
$^3$ Кафедра прикладной математики и программирования, Южно-Уральский государственный университет, Челябинск, Россия\\
$^4$ Институт систем энергетики им. Л.А. Мелентьева, Сибирское отделение
Российской академии наук,  Иркутск, Россия\\
$^5$ Институт математики и информационных технологий, Иркутский государ-ственный университет, Иркутск, Россия
}}

%аннотация,
\begin{abstracte} Предложен метод полиномиальной сплайн-коллокации для решения интегральных уравнений Вольтерра первого рода с  кусочно-непрерывными ядрами.
Для аппроксимации интегралов при дискретизации в предлагаемом проекционном методе используется квадратурная формула типа Гаусса. Получена оценка точности приближенного решения. Также используется стохастическая арифметика (СА) на основе метода Controle et Estimation Stochastique des Arrondis de Calculs (CESTAC) и библиотеки Control of Accuracy and Debugging for Numerical Applications (CADNA). Применяя этот подход, можно найти оптимальные параметры проекционного метода. Приведены численные примеры, иллюстрирующие эффективность предлож\-енного нового метода коллокации.
\end{abstracte}

%ключевые слова
\keywordse{интегральное уравнение, разрывное ядро, метод сплайн-коллокации, сходимость, метод CESTAC, библиотека CADNA.}
\selectlanguage{russian}
%список литературы

%
%%Сведения об авторах
%\textbf{Иван Иванович Иванов}, доктор физико-математических наук, профессор, Иркутский государственный университет, Российская федерация, Иркутск, 664003, ул. К. Маркса, 1,
%тел.: +7(3952)242210,\\
%email: avtor@math.isu.ru, ORCID iD https://orcid.org/xxxx-xxxx-xxxx-xxxx.
%%\selectlanguage{russian}
%
%\textbf{Петр Петрович Петров}, кандидат физико-математических наук,
%Новосибирский государственный университет,
%Российская федерация,  Новосибирск, 630090, ул. Пирогова, 1,
%тел.:+7 (383)3634333,\\
%email: avtor@math.isu.ru, ORCID iD https://orcid.org/xxxx-xxxx-xxxx-xxxx.
%

\textbf{Александр Николаевич Тында}, кандидат физико-математических наук, доцент, заместитель декана по учебной работе, кафедра высшей и прикладной математики, Пензенский государственный университет, 440026 Пенза, Россия; \\
email: tyndaan@mail.ru, \\
ORCID ID https://orcid.org/0000-0001-6023-9847.
\vspace{0.3cm}

\textbf {Самад Нойягдам}, доктор (PhD) прикладной математики, доцент,
Лаборатория промышленной математики Байкальской школы БРИКС Иркутского национального исследовательского технического университета, Иркутск, 664074, Россия; \\
Кафедра прикладной математики и программирования, Южно-Уральский государственный университет, Россия, 454080, г. Челябинск, проспект Ленина, 76; \\
 email: snoei@istu.edu; noiagdams@susu.ru; \\
 ORCID ID https://orcid.org/0000-0002-2307-0891.
\vspace{0.3cm}

\textbf{Денис Николаевич Сидоров}, доктор физико-математических наук, профессор РАН, главный научный сотрудник,  Институт  систем энергетики им. Л.А. Мелентьева СО РАН, 664033 Иркутск, Россия; \\
Зав. лаб. промышленной математики Байкальской школы БРИКС Иркутского национального исследовательского технического университета, Иркутск, 664074, Россия; \\
Институт математики и информационных технологий Иркутского государственного университета, Иркутск 664025, Россия;
 email: contact.dns@gmail.com, \\ ORCID ID https://orcid.org/0000-0002-3131-1325.

\end{article}

\end{document}